# Solving the OSCAR and SLOPE Models Using a Semismooth Newton-Based Augmented Lagrangian Method

**Ziyan Luo**                                              STARKEYNATURE@HOTMAIL.COM
*State Key Laboratory of Rail Traffic Control and Safety*
*Beijing Jiaotong University*
*Beijing, P. R. China*

**Defeng Sun**                                            DEFENG.SUN@POLYU.EDU.HK
*Department of Applied Mathematics*
*The Hong Kong Polytechnic University*
*Hong Kong*

**Kim-Chuan Toh**                                         MATTOHKC@NUS.EDU.SG
*Department of Mathematics & Institute of Operations Research and Analytics*
*National University of Singapore*
*Singapore*

**Naihua Xiu**                                            NHXIU@BJTU.EDU.CN
*Department of Mathematics*
*Beijing Jiaotong University*
*Beijing, P. R. China*

## Abstract

The octagonal shrinkage and clustering algorithm for regression (OSCAR), equipped with the $\ell_1$-norm and a pair-wise $\ell_\infty$-norm regularizer, is a useful tool for feature selection and grouping in high-dimensional data analysis. The computational challenge posed by OSCAR, for high dimensional and/or large sample size data, has not yet been well resolved due to the non-smoothness and inseparability of the regularizer involved. In this paper, we successfully resolve this numerical challenge by proposing a sparse semismooth Newton-based augmented Lagrangian method to solve the more general SLOPE (the sorted L-one penalized estimation) model. By appropriately exploiting the inherent sparse and low-rank property of the generalized Jacobian of the semismooth Newton system in the augmented Lagrangian subproblem, we show how the computational complexity can be substantially reduced. Our algorithm presents a notable advantage in the high-dimensional statistical regression settings. Numerical experiments are conducted on real data sets, and the results demonstrate that our algorithm is far superior, in both speed and robustness, than the existing state-of-the-art algorithms based on first-order iterative schemes, including the widely used accelerated proximal gradient (APG) method and the alternating direction method of multipliers (ADMM).

**Keywords:** Linear Regression, OSCAR, Sparsity, Augmented Lagrangian Method, Semismooth Newton method

# 1. Introduction

Feature selection and grouping is highly beneficial in learning with high-dimensional data containing spurious features, and thus has found wide applications in computer vision, signal processing, bioinformatics, etc. The octagonal shrinkage and clustering algorithm for regression (OSCAR) proposed by Bondell and Reich (2008), serves as an efficient sparse modeling tool with automatic feature grouping by employing the $\ell_1$-norm regularizer together with a pairwise $\ell_\infty$ penalty. The OSCAR penalized problem for linear regression with the least squares loss function takes the form of

$$\min_{x \in \mathbb{R}^n} \frac{1}{2}\|Ax - b\|_2^2 + w_1\|x\|_1 + w_2 \sum_{i<j} \max\left\{|x_i|, |x_j|\right\}, \tag{1}$$

where $b \in \mathbb{R}^m$ is the response vector, $A \in \mathbb{R}^{m \times n}$ is the design matrix, $x \in \mathbb{R}^n$ is the vector of unknown coefficients to be estimated, $w_1$ and $w_2$ are two nonnegative tuning parameters for the tradeoff of the sparsity and equality of coefficients for correlated features promoted by the $\ell_1$-norm and the pairwise $\ell_\infty$ term, respectively. In the high dimensional setting of statistical regression, we always have $n > m$, that is, the number of features is larger than the sample size.

The OSCAR penalized problem (1) is a convex optimization problem. When the pairwise $\ell_\infty$ term is removed, the problem (1) is reduced to the well-known LASSO model proposed by Tibshirani (1996) in statistics and a rich variety of algorithms have been proposed, most of which have taken the advantage of the componentwise separability of the $\ell_1$-norm in their algorithmic design. With the additional pairwise $\ell_\infty$ term, the problem (1) becomes understandably more challenging due to the lack of separability of the OSCAR regularization term. In Bondell and Reich (2008), the traditional quadratic programming (QP) and sequential quadratic programming (SQP) based algorithms are employed for solving (1) with numerical implementations limited to small data sets . Efficient numerical algorithms are in dire need especially for large scale problems resulting from the explosion in the size and complexity of modern data sets in practical applications. In Zhong and Kwok (2012), the accelerated proximal gradient (APG) method, proposed by Nesterov (1983) and coined as FISTA for the $\ell_1$-norm regularization problem by Beck and Teboulle (2009), is adopted for solving relatively large scale instances by taking advantage of the efficient computation of the proximal mapping of the OSCAR penalty function. Note that the OSCAR penalty can be written as

$$w_1\|x\|_1 + w_2 \max_{i<j}\left\{|x_i|, |x_j|\right\} = \sum_{i=1}^n \lambda_i |x|_i^\downarrow \tag{2}$$

by using the non-increasing components in magnitude $|x|_1^\downarrow \geq \cdots \geq |x|_n^\downarrow$ and the non-increasing parameters $\lambda_i = w_1 + w_2(n-i)$, $i = 1, \ldots, n$. The resulting regularization function $\kappa_\lambda(x) := \sum_{i=1}^n \lambda_i |x|_i^\downarrow$ with $\lambda_1 \geq \lambda_2 \cdots \geq \lambda_n \geq 0$ for any $x \in \mathbb{R}^n$, termed as the decreasing weighted sorted $\ell_1$ norm (DWSL1) by Zeng and Figueiredo (2014), is exactly the weighted Ky Fan norm as studied in Wu et al. (2014) as long as $\lambda_1 > 0$. The computation of the proximal mapping of DWSL1 has been studied in the literature (see, e.g., Zeng and Figueiredo, 2013, 2014; Bogdan et al., 2015). It is heavily related to the pool adjacent violators algorithm (PAVA) for solving isotonic regression problems (Barlow and Brunk, 1972) in the field of ordered statistics (see, e.g., Robertson et al., 1988; Silvapulle and Sen, 2011).

As a more general framework of the OSCAR problem (1), the least-squares problem with the DWSL1 regularization term is called the sorted L-one penalized estimation (SLOPE), which has been shown to have nice performance for controlling the false discovery rate (FDR) in sparse statistical models as pointed in Bogdan et al. (2015). The APG method is then utilized for solving the SLOPE model by relying on the efficient numerical evaluation of the proximal mapping of the involved sorted $\ell_1$ norm. As can be seen, most of the existing methods for solving the OSCAR and its general case SLOPE in the large scale settings are based on the first-order information of the underlying nonsmooth optimization model. However, as demonstrated by the works of Li et al. (2018) for the LASSO and Li et al. (2017) for the fused LASSO, there are compelling evidences to suggest that one can design a much more efficient algorithm if one can fully exploit the inherent second-order sparsity and low rank property present in the OSCAR model or more generally the SLOPE model. In this paper, we will show how this can be achieved by focusing on the following SLOPE model

$$\min_{x \in \mathbb{R}^n} \frac{1}{2} \|Ax - b\|_2^2 + \sum_{i=1}^n \lambda_i |x|_i^{\downarrow} \tag{3}$$

with parameters $\lambda_1 \geq \lambda_2 \geq \cdots \geq \lambda_n \geq 0$ and $\lambda_1 > 0$. A semismooth Newton-based augmented Lagrangian method (Newt-ALM for short) will be applied for solving the SLOPE model from the dual perspective. By wisely taking advantage of the special structure of the Hessian matrices involved in the inner iterations and the fast linear convergence of the augmented Lagrangian method for the outer iterations, the Newt-ALM is demonstrated to perform highly efficient in numerical experiments on large scale instances. The comparison of our algorithm with the inexact ADMM (iADMM) proposed in Chen et al. (2017) and the SLOPE solver in Bogdan et al. (2015) for solving OSCAR problems indicates that our Newt-ALM outperforms other existing algorithms substantially. A key contribution to the high computational efficiency of our proposed Newt-ALM is the characterization of the generalized Jacobian matrix of the involved Newton system in each subproblem, where the inherent low-rank and sparsity structure resulting from the proximal mapping of the DWSL1 term can be fully exploited to greatly reduce the computational cost. This constitutes one of the main contributions of this paper.

The remaining parts of the paper are organized as follows. In Section 2, selected properties on proximal mappings and their generalized Jacobians are reviewed and developed, which are critical for the subsequent analysis on local convergence rates. Section 3 is dedicated to the semismooth Newton augmented Lagrangian method and its convergence analysis. Numerical results are reported in Section 4 to show the high efficiency and robustness of our algorithm, and technical proofs are provided in Appendix A. We conclude our paper in Section 5.

## 2. The generalized Jacobian of the proximal mapping of the DWSL1 function

Let $\mathbf{\Pi_n^s}$ be the set of all signed permutation matrices in $\mathbb{R}^{n \times n}$. Obviously, the cardinality of $\mathbf{\Pi_n^s}$ is $2^n n!$. For any given vector $y \in \mathbb{R}^n$, denote

$$\Pi^s(y) := \left\{ \pi \in \mathbf{\Pi_n^s} \mid \pi(\mathbf{y}) = |\mathbf{y}|^\downarrow \right\},$$

where $|z|^\downarrow$ stands for the vector of components in $|z|$ being arranged in the non-increasing order. Let $\kappa_\lambda(x) := \sum_{i=1}^n \lambda_i |x|_i^\downarrow$ with $\lambda_1 \geq \cdots \geq \lambda_n \geq 0$. The proximal mapping of $\kappa_\lambda$ is

$$\mathrm{Prox}_{\kappa_\lambda}(y) = \arg\min_x \left\{ \frac{1}{2} \|x - y\|^2 + \kappa_\lambda(x) \right\}, \quad \forall y \in \mathbb{R}^n.$$

Since the involved objective function is strongly convex (see, e.g., Wu et al., 2014; Bogdan et al., 2015) and piecewise quadratic, the proximal mapping $\mathrm{Prox}_{\kappa_\lambda}$ is then piecewise affine, a result known from Sun (1986) or (Rockafellar and Wets, 1998, Proposition 12.30). Define

$$x_\lambda(w) := \arg\min_x \left\{ \frac{1}{2} \|x - w\|^2 + \lambda^\top x \mid Bx \geq 0 \right\}, \quad w \in \Re^n,$$

where

$$Bx = [x_1 - x_2, x_2 - x_3, \ldots, x_{n-1} - x_n, x_n]^\top \in \mathbb{R}^n.$$

It is known from (Bogdan et al., 2015, Proposition 2.2) that for any $y \in \mathbb{R}^n$ and $\pi \in \Pi^s(y)$, $\mathrm{Prox}_{\kappa_\lambda}(\pi y) = x_\lambda(\pi y)$. Furthermore, for any $\lambda \in \mathbb{R}_+^n$ satisfying $\lambda = |\lambda|^\downarrow$, and any vector $y \in \mathbb{R}^n$, we have

$$\mathrm{Prox}_{\kappa_\lambda}(y) = \pi^{-1} x_\lambda(\pi y), \quad \forall \pi \in \Pi^s(y) \subseteq \mathbf{\Pi_n^s}. \tag{4}$$

Given the structure of $x_\lambda(\cdot)$, one can see that the HS-Jacobian of $x_\lambda(\cdot)$ at any $w \in \mathbb{R}^n$, (see, e.g., Han and Sun, 1997), termed as $\mathcal{P}(w)$, is defined by

$$\mathcal{P}(w) = \left\{ P \in \mathbb{R}^{n \times n} \mid P = I - B_\Gamma^\top \left( B_\Gamma B_\Gamma^\top \right)^{-1} B_\Gamma, \Gamma \in \mathcal{K}(w) \right\}, \tag{5}$$

where $\mathcal{K}(w) := \{ \Gamma \subseteq \{1, \ldots, n\} \mid \mathrm{Supp}(z_\lambda(w)) \subseteq \Gamma \subseteq \mathbf{I}(x_\lambda(w)) \}$ with the optimal dual solution $z_\lambda(w) = (BB^\top)^{-1} B(w - \lambda - x_\lambda(w))$ and $\mathbf{I}(x_\lambda(w)) = \{ i \in \{1, \ldots, n\} \mid (Bx_\lambda(w))_i = 0 \}$ and $B_\Gamma$ is the submatrix obtained by extracting the rows of $B$ with indices in $\Gamma$. Here, $\mathrm{Supp}(z_\lambda(w))$ is the support of $z_\lambda(w)$, i.e., the index set of nonzero compomonents of $z_\lambda(w)$.

It is known from Lemma 2.1 in Han and Sun (1997) that for any $w \in \mathbb{R}^n$, there exists a neighborhood $W$ of $w$ such that for all $w' \in W$,

$$\begin{cases} \mathcal{K}(w') \subseteq \mathcal{K}(w), \\ \mathcal{P}(w') \subseteq \mathcal{P}(w), \\ x_\lambda(w') - x_\lambda(w) - P(w' - w) = 0, \quad \forall P \in \mathcal{P}(w'). \end{cases} \tag{6}$$

Define the multifunction $\mathcal{M} : \mathbb{R}^n \rightrightarrows \mathbb{R}^{n \times n}$ by

$$\mathcal{M}(y) := \left\{ M \in \mathbb{R}^{n \times n} \mid M = \pi^{-1} P \pi, \pi \in \Pi^s(y), P \in \mathcal{P}(\pi y) \right\}. \tag{7}$$

Then we have the following theorem, whose proof is given in Appendix A.

**Theorem 1** *Let $\lambda \in \mathbb{R}^n_+$ be such that $\lambda = |\lambda|^{\downarrow}$. Then $\mathcal{M}(\cdot)$ is a nonempty and compact valued, upper semicontinuous multifunction, and for any given $y \in \mathbb{R}^n$, every $M \in \mathcal{M}(y)$ is symmetric and positive semidefinite. Moreover, there exists a neighborhood $U$ of $y$ such that for all $y' \in U$,*

$$\mathrm{Prox}_{\kappa_\lambda}(y') - \mathrm{Prox}_{\kappa_\lambda}(y) - M(y' - y) = 0, \quad \forall M \in \mathcal{M}(y'). \tag{8}$$

Recall from Mifflin (1977); Kummer (1988); Qi and Sun (1993); Sun and Sun (2002) or directly from (Li et al., 2017, Definition 1) that the semismoothness with respect to a given nonempty compact-valued, upper semicontinuous multifunction is defined as follows.

Let $\mathcal{O} \subseteq \mathbb{R}^n$ be any given open set, $\mathcal{K} : \mathcal{O} \rightrightarrows \mathbb{R}^{m \times n}$ be a nonempty compact valued, upper semicontinuous multifunction, and $F : \mathcal{O} \to \mathbb{R}^m$ be a locally Lipschitz continuous function. $F$ is said to be *semismooth* at $x \in \mathcal{O}$ with respect to the multifunction $\mathcal{K}$ if $F$ is directionally differentiable at $x$ and for any $V \in \mathcal{K}(x + d)$ with $d \to 0$,

$$F(x + d) - F(x) - Vd = o(\|d\|).$$

Let $\gamma$ be a positive scalar. $F$ is said to be $\gamma$-order semismooth (stongly semismooth if $\gamma = 1$) at $x \in \mathcal{O}$ with respect to $\mathcal{K}$ if $F$ is directionally differentiable at $x$ and for any $V \in \mathcal{K}(x + d)$ with $d \to 0$,

$$F(x + d) - F(x) - Vd = O(\|d\|^{1+\gamma}).$$

$F$ is said to be a semismooth ($\gamma$-order semismooth, stongly) function on $\mathcal{O}$ with respect to $\mathcal{K}$ if $F$ is semismooth ($\gamma$-order semismooth, strongly semismooth) everywhere in $\mathcal{O}$ with respect to $\mathcal{K}$. It is known from Theorem 1 that $\mathrm{Prox}_{\kappa_\lambda}$ is $\gamma$-order semismooth on $\mathbb{R}^n$ with respect to $\mathcal{M}$ for any given positive $\gamma$.

## 3. A semismooth Newton augmented Lagrangian method

### 3.1 The algorithmic framework

Given $A \in \mathbb{R}^{m \times n}$, $b \in \mathbb{R}^m$ and $\lambda_1 \geq \cdots \geq \lambda_n \geq 0$ and $\lambda_1 > 0$, the DWSL1 regularized least squares problem can be rewritten as

$$(P) \quad \max_{x \in \mathbb{R}^n} \left\{ -f(x) := -\frac{1}{2}\|Ax - b\|^2 - \kappa_\lambda(x) \right\}. \tag{9}$$

Its dual problem takes the form of

$$(D) \quad \min_{y \in \mathbb{R}^m} \left\{ g(y) := \frac{1}{2}\|y\|^2 + \langle b, y \rangle + \kappa_\lambda^*(-A^\top y) \right\}, \tag{10}$$

where $\kappa_\lambda^*(v) := \sup_{x \in \mathbb{R}^n} \{\langle x, v \rangle - \kappa_\lambda(x)\}$ is the Fenchel conjugate function of $\kappa_\lambda$. Following from the same scheme in (Rockafellar and Wets, 1998, Examples 11.46 and 11.57) (also Section 4 Li et al., 2017), the Lagrangian function $l : \mathbb{R}^n \times \mathbb{R}^n \to [-\infty, +\infty]$ associated with $(D)$ is given by

$$\begin{aligned} l(y; x) &:= \inf_u \left\{ \frac{1}{2}\|y\|^2 + \langle b, y \rangle + \kappa_\lambda^*(-A^\top y + u) \right\} \\ &= \frac{1}{2}\|y\|^2 + \langle b, y \rangle - \langle A^\top y, x \rangle - \kappa_\lambda(x). \end{aligned} \tag{11}$$

For any given scalar $\sigma > 0$, the corresponding augmented Lagrangian function associated with $(D)$ is defined by

$$
\begin{aligned}
L_\sigma(y; x) &:= \inf_u \left\{ \frac{1}{2} \|y\|^2 + \langle b, y \rangle + \kappa_\lambda^*(-A^\top y + u) + \frac{\sigma}{2} \|u\|^2 \right\} \\
&= \frac{1}{2} \|y\|^2 + \langle b, y \rangle + \inf_{s \in \mathbb{R}^n} \left\{ \kappa_\lambda^*(s) - \langle x, A^\top y + s \rangle + \frac{\sigma}{2} \|A^\top y + s\|^2 \right\} \\
&= \frac{1}{2} \|y\|^2 + \langle b, y \rangle - \frac{1}{2\sigma} \|x\|^2 + \sigma \phi_{\kappa_\lambda^*/\sigma} \left( \frac{x + \sigma A^\top y}{\sigma} \right),
\end{aligned}
\tag{12}
$$

where $\phi_{\kappa_\lambda^*/\sigma}$ is the Moreau-Yosida regularization of $\kappa_\lambda^*/\sigma$ defined as

$$
\phi_{\kappa_\lambda^*/\sigma}(x) := \min_{u \in \mathbb{R}^n} \left\{ \frac{1}{\sigma} \kappa_\lambda^*(u) + \frac{1}{2} \|u - x\|^2 \right\}, \quad \forall x \in \mathbb{R}^n.
$$

The inexact augmented Lagrangian method (Rockafellar, 1976b) together with the semismooth Newton method will be employed to solve $(D)$ with the algorithmic framework as described in Algorithm 1.

---

**Algorithm 1:** An inexact augmented Lagrangian method for $(D)$ (Newt-ALM)

Choose $\sigma_0 > 0$ and $(y^0, x^0) \in \mathbb{R}^m \times \mathbb{R}^n$. For $k = 0, 1, \ldots,$ perform the following steps in each iteration:

**Step 1.** Compute $y^{k+1} \approx \arg \min_{y \in \mathbb{R}^m} \left\{ \Psi_k(y) := L_{\sigma_k}(y; x^k) \right\}$;

**Step 2.** Compute $x^{k+1} = \mathrm{Prox}_{\sigma_k \kappa_\lambda} \left( x^k - \sigma_k A^\top y^{k+1} \right)$;

**Step 3.** Update $\sigma_{k+1} \uparrow \sigma_\infty \leq \infty$.

---

The stopping criteria for the inexact augmented Lagrangian method have been well discussed in Rockafellar (1976b,a). Given two summable sequences of nonnegative numbers, $\{\epsilon_k\}_{k \geq 0}$ and $\{\delta_k\}_{k \geq 0}$, and a nonnegative convergent sequence $\{\delta_k'\}_{k \geq 0}$ with limit 0, the stopping criteria can be simplified as follows in our case:

(A) $\|\nabla \Psi_k(y^{k+1})\| \leq \epsilon_k / \sqrt{\sigma_k}$;

(B1) $\|\nabla \Psi_k(y^{k+1})\| \leq (\delta_k / \sqrt{\sigma_k}) \|x^{k+1} - x^k\|$;

(B2) $\|\nabla \Psi_k(y^{k+1})\| \leq (\delta_k' / \sigma_k) \|x^{k+1} - x^k\|$.

## 3.2 Convergence Theory

The piecewise linear-quadratic property of $f$ as defined in (9) leads to the polyhedral multifunction $\partial f$ (the sub-differential of $f$), which further implies that $\partial f$ satisfies the error bound condition with a common modulus, say $a_f$. Especially, since the optimal solution set of (P), denoted by $S^*$, is exactly $(\partial f)^{-1}(0)$, there exists some $\varepsilon > 0$ such that for any $x \in \mathbb{R}^n$ satisfying $\mathrm{dist}(0, \partial f(x)) \leq \varepsilon$, it holds that

$$
\mathrm{dist}(x, S^*) \leq a_f \, \mathrm{dist}(0, \partial f(x)).
\tag{13}
$$

Similarly, for the polyhedral multifunction $T_l$ defined as $T_l(y, x) = \{(y', x') \mid (y', -x') \in \partial l(y; x)\}$, there exist some $a_l$ and $\varepsilon' > 0$ such that for any $(y, x) \in \mathbb{R}^m \times \mathbb{R}^n$ satisfying $\text{dist}(0, T_l(y, x)) \leq \varepsilon'$, it has

$$\text{dist}\left((y, x), \{y^*\} \times S^*\right) \leq a_l \, \text{dist}(0, T_l(y, x)), \tag{14}$$

where $y^*$ is the unique optimal solution of (D). Following the results on global and local convergence of the ALM as stated in (Rockafellar, 1976b,a; Li et al., 2018, 2017), we can readily obtain the following convergence results on Algorithm 1 with the above criteria.

**Theorem 2 (Global convergence)** *Let $\left\{(y^k, x^k)\right\}$ be the infinite sequence generated by Algorithm 1 with stopping criterion (A). Then $\left\{x^k\right\}$ converges to an optimal solution to (P), and $\left\{y^k\right\}$ converges to the unique optimal solution of (D).*

**Theorem 3 (Local linear-rate convergence)** *Let $\left\{(y^k, x^k)\right\}$ be the infinite sequence generated by Algorithm 1 with stopping criteria (A) and (B1). Then for all $k$ sufficiently large,*

$$\text{dist}\left(x^{k+1}, S^*\right) \leq \theta_k \text{dist}\left(x^k, S^*\right),$$

*where*

$$\theta_k = \left(\frac{a_f}{\sqrt{a_f^2 + \sigma_k^2} + 2\delta_k}\right) / (1 + \delta_k) \to \theta_\infty := \frac{a_f}{\sqrt{a_f^2 + \sigma_\infty^2}} < 1$$

*as $k \to +\infty$, and $a_f$ is from (13). Additionally, if the criterion (B2) is also adopted, then for all $k$ sufficiently large,*

$$\|y^{k+1} - y^*\| \leq \theta' \|x^{k+1} - x^k\|,$$

*where*

$$\theta_k' = \frac{a_l(1 + \delta_k')}{\sigma_k} \to \frac{a_l}{\sigma_\infty}$$

*as $k \to +\infty$, and $a_l$ is from (14).*

**Remark 4 (Global linear-rate convergence)** *Besides the local linear-rate convergence as stated in Theorem 3, one can also obtain the global Q-linear convergence of the primal sequence $\{x^k\}$ and the global R-linear convergence of the dual infeasibility and the duality gaps for the sequence generated by Algorithm 1 based on (Cui et al., 2017, Proposition 2 and Lemma 3) or by mimicking the proofs of (Zhang et al., 2017, Theorem 4.1 and Remark 4.1) since problem (P) possesses the following property: For any positive scalar $r$, there exists $t > 0$ such that*

$$\text{dist}(x, S^*) \leq t \, \text{dist}(0, \partial f(x)), \;\; \forall x \in \mathbb{R}^n \text{ satisfying } \text{dist}(x, S^*) \leq r, \tag{15}$$

*(see, Zhang et al., 2017, Proposition 2.2). We omit the details here.*

### 3.3 The semismooth Newton method for the essential subproblem

It is known from Moreau (1965) or (Rockafellar, 1970, Theorem 31.5) that $\Psi_k$ is continuously differentiable and

$$\nabla\Psi_k(y) = y + b - A\mathrm{Prox}_{\sigma_k\kappa_\lambda}(x^k - \sigma_k A^\top y), \quad \forall y \in \mathbb{R}^m. \tag{16}$$

Since $\Psi_k$ is strongly convex with bounded level sets, the unique solution of $\min_{y\in\mathbb{R}^m}\Psi_k(y)$ can be computed by the following first-order optimality condition

$$\nabla\Psi_k(y) = 0. \tag{17}$$

For any $y \in \mathbb{R}^n$, define

$$\mathcal{G}_k(y) := \left\{ V \in \mathbb{R}^{m\times m} \mid V = I_m + \sigma_k AMA^\top, \; M \in \mathcal{M}(x^k - \sigma_k A^\top y) \right\},$$

where $\mathcal{M}$ is defined in (7). The following semismooth Newton (SSN) method is then applied to solve the semismooth equation (17), as presented in Algorithm 2.

---

**Algorithm 2:** A semismooth Newton method for solving (17)

---

Choose $\mu \in (0, 1/2)$, $\bar{\eta} \in (0, 1)$, $\tau \in (0, 1]$, $y^0 \in \mathbb{R}^m$. For $j = 0, 1, \ldots$, perform the following steps in each iteration:

**Step 1.** (Computing the Newton direction) Let $M_j$ be an element in $\mathcal{M}(x^k - \sigma_k A^\top y^j)$ and set $V_j := I_m + \sigma_k AM_jA^\top$. Solve the Newton equation

$$V_j d = -\nabla\Psi_k(y^j) \tag{18}$$

exactly or by the conjugate gradient (CG) algorithm to get $d^j$ such that $\|V_j d^j + \nabla\Psi(y^j)\| \le \min\{\bar{\eta}, \|\nabla\Psi(y^j)\|^{1+\tau}\}$.

**Step 2.** (Line search) Set $\alpha_j = \delta^{m_j}$, where $m_j$ is the least nonnegative integer $m$ satisfying

$$\Psi(y^j + \delta^m d^j) \le \Psi(y^j) + \mu\delta^m\langle\nabla\Psi(y^j), d^j\rangle.$$

**Step 3.** Set $y^{j+1} = y^j + \alpha_j d^j$.

---

### 3.4 Efficient implementations of the semismooth Newton method

In this subsection, the sparsity and low-rank structure of the coefficient matrix in the linear system (18) will first be uncovered. Then the structures will be exploited through designing novel numerical techniques for solving the large scale system (18) to achieve efficient implementations of the semismooth Newton method in Algorithm 2. For any given index set $\Gamma \subseteq \{1, \ldots, n\}$, define the diagonal matrix $\Sigma_\Gamma \in \mathbb{R}^{n\times n}$ by

$$(\Sigma_\Gamma)_{ii} = \begin{cases} 1, & \text{if } i \in \Gamma; \\ 0, & \text{otherwise.} \end{cases}$$

Similar to the case in (Li et al., 2017, Proposition 6), there exists some positive integer $N$ such that $\Sigma_\Gamma$ can be rewritten as a block diagonal matrix

$$\Sigma_\Gamma = \mathrm{Diag}(\Lambda_1, \ldots, \Lambda_N)$$

with $\Lambda_i \in \{O_{n_i}, I_{n_i}\}$ for each $i \in \{1, \ldots, N\}$ where any two consecutive blocks $\Lambda_i$ and $\Lambda_{i+1}$ are not of the same type. Denote

$$J = \{j \in \{1, \ldots, N\} \mid \Lambda_j = I\}.$$

Then we have

$$P = I_n - B_\Gamma^\top (B_\Gamma B_\Gamma^\top)^{-1} B_\Gamma = \mathrm{Diag}(P_1, \ldots, P_N),$$

where

$$P_i = \begin{cases} \frac{1}{n_i+1} e_{n_i+1} e_{n_i+1}^\top, & \text{if } i \in J \text{ and } i \neq N; \\ O_{n_i}, & \text{if } i \in J \text{ and } i = N; \\ I_{n_i-1}, & \text{if } i \notin J \text{ and } i \neq 1; \\ I_{n_i}, & \text{if } i \notin J \text{ and } i = 1 \end{cases}$$

with the convention $I_0 = \emptyset$. This block diagonal matrix $P$ can be further decomposed into the sum of a sparse diagonal term and a low-rank term as $P = H + UU^T$, where $H = \mathrm{Diag}(H_1, \ldots, H_N) \in \mathbb{R}^{n \times n}$ with

$$H_i = \begin{cases} O_{n_i+1}, & \text{if } i \in J \text{ and } i \neq N; \\ O_{n_i}, & \text{if } i \in J \text{ and } i = N; \\ I_{n_i-1}, & \text{if } i \notin J \text{ and } i \neq 1; \\ I_{n_i}, & \text{if } i \notin J \text{ and } i = 1 \end{cases}$$

and $U \in \mathbb{R}^{n \times N}$ with its $(k, j)^{\text{th}}$ entry

$$U_{kj} = \begin{cases} 1/\sqrt{n_j + 1}, & \text{if } \sum\limits_{t=1}^{j-1} n_t + 1 \leq k \leq \sum\limits_{t=1}^{j} n_t + 1 \text{ and } j \in J \backslash \{N\}; \\ 0, & \text{otherwise.} \end{cases}$$

Define $\alpha := \{j \in \{1, \ldots, n\} \mid H_{ii} = 1\} = \{1, \ldots, n\} \backslash \Gamma$, and denote $U_{J_N}$ as the submatrix of $U$ generated from columns of $U$ indexed by $J \backslash \{N\}$. Then for any given $A \in \mathbb{R}^{m \times n}$, any $\Gamma \in \{1, \ldots, n\}$ with its corresponding matrix $P$ defined as above, and any signed permutation matrix $\pi$, we have

$$\begin{aligned} A\pi^\top P \pi A^\top &= A\pi^\top H \pi A^\top + A\pi^\top U U^\top \pi A^\top \\ &= A\pi(\alpha, :)^\top \pi(\alpha, :) A^\top + \widetilde{A} \widetilde{U}_{J_N} \widetilde{U}_{J_N}^\top \widetilde{A}^\top \\ &=: V_1 V_1^\top + V_2 V_2^\top, \end{aligned} \tag{19}$$

where $V_1 = A\pi(\alpha, :)^\top$, $V_2 = \widetilde{A} \widetilde{U}_{J_N}$ with $\widetilde{U}_{J_N}$ being the submatrix of $U_{J_N}$ obtained by dropping all its zero rows and $\widetilde{A}$ is the submatrix obtained from the permuted matrix $A\pi^\top$ by dropping the columns corresponding to those zero rows in $U_{J_N}$. Therefore, the cost of computing $A\pi^{-1} P \pi A^\top$ is dramatically reduced from $\mathcal{O}(mn(n+m))$ by naive computation to $\mathcal{O}(m^2(r_1 + r_2))$, where $r_1$ is number of columns in $V_1$ and $r_2$ is the number of columns in $V_2$. Here $r_1$ refers to the number of inactive constraints in $Bx \geq 0$, and $r_2$ refers to the number of distinct nonzero identical components in $Bx$, both of which are generally no larger than the number of nonzero components of $x$. In the setting of high-dimensional sparse grouping linear regression models, $m$, $r_1$, $r_2$ and $N$ are generally much smaller than $n$, therefore the aforementioned reduction of the computational cost can be highly significant.

If $m$ is not too large, we can use the (sparse) Cholesky factorization to directly solve the linear system (17). In the case where $r_1 + r_2 \ll m$, the cost of solving (17) can be further reduced by using the Sherman-Morrison-Woodbury formula as follows:

$$\left( I_m + \sigma A \pi^{-1} P \pi A^\top \right)^{-1} = (I_m + WW^\top)^{-1} = I_m - W(I_{r_1+r_2} + W^\top W)^{-1} W^\top,$$

where $W = \sqrt{\sigma}[V_1 \quad V_2] \in \mathbb{R}^{m \times (r_1+r_2)}$. In the event when $m$ is extremely large, we can use the preconditioned conjugate gradient (PCG) method to solve the linear system (17).

## 4. Numerical experiments

The performance of our proposed algorithm for solving SLOPE (3) and the special case of the OSCAR model in (1) will be evaluated by comparing with the accelerated proximal gradient (APG) algorithm implemented in Bogdan et al. (2015) with its Matlab code available at http://statweb.stanford.edu/ candes/SortedL1, the semi-proximal alternating direction method of multipliers (sPADMM) (see, e.g., Fazel et al. (2013)) applied to the dual problem (with implementation details presented in Subsection 4.2), and our sparse semismooth Newton-based augmented Lagrangian method (Newt-ALM). All the computational results are obtained from a desktop computer running on 64-bit Windows Operating System having 4 cores with Intel(R) Core(TM) i5-5257U CPU at 2.70GHz and 8 GB memory.

### 4.1 Stopping Criteria

To measure the accuracy of an approximate optimal solution $(y, x)$ for the dual problem (10) and the primal problem (9), the relative duality gap and the dual infeasibility will be adopted. Specifically, denote

$$Obj_P := \frac{1}{2}\|Ax - b\|^2 + \kappa_\lambda(x), \text{ and } Obj_D := -b^\top y - \frac{1}{2}\|y\|^2,$$

which are the values of the primal and dual objective functions, respectively. Then the relative duality gap can be defined by

$$\eta_G := \frac{|Obj_P - Obj_D|}{\max\{1, |Obj_P|\}}.$$

Note that $\kappa_\lambda^*(\cdot)$ is actually the indicator function induced by the closed convex set

$$\mathcal{C}_\lambda := \left\{ z \ \Big| \ \sum_{j \leq i} |z|_j^\downarrow \leq \sum_{j \leq i} \lambda_j, \ i = 1, \dots, n \right\},$$

which is exactly the unit ball of the dual norm to $\kappa_\lambda$ (see, e.g., Bogdan et al., 2015; Wu et al., 2014). To characterize the dual infeasibility of $y$, or equivalently $-A^\top y \in \mathcal{C}_\lambda$, we adopt the measurement proposed in Bogdan et al. (2015) which is

$$\eta_D := \max \left\{ 0, \max_{1 \leq i \leq n} \sum_{j \leq i} \left( |A^\top(Ax - b)|_j^\downarrow - \lambda_j \right) \right\}.$$

For given accuracy parameters $\varepsilon_G$ and $\varepsilon_D$, our algorithm Newt-ALM will be terminated once

$$\eta_G \leq \varepsilon_G \text{ and } \eta_D \leq \varepsilon_D, \tag{20}$$

while both the sPADMM and the APG will be terminated if (20) holds or if the number of iterations reaches the maximum of $50,000$. In our numerical experiments, we choose $\varepsilon_G = \varepsilon_D = $ 1e-6.

The relative KKT residual

$$\eta = \frac{\|x - \text{Prox}_{\kappa_\lambda}(x - A^\top(Ax - b))\|}{1 + \|x\| + \|A^\top(Ax - b)\|}$$

is adopted to measure the accuracy of an approximate optimal solution $x$ from any of the algorithm implemented in the numerical experiments.

### 4.2 ADMM for problem (10)

The implementation details of the (semi-proximal) ADMM for solving problem (10) are elaborated in this subsection. By reformulating (10) as the following equality-constrained problem

$$(D') \quad \min_{y \in \mathbb{R}^m} \left\{ \frac{1}{2}\|y\|^2 + \langle b, y \rangle + \kappa_\lambda^*(\xi) \, \Big| \, A^\top y + \xi = 0 \right\}, \tag{21}$$

we can state a general framework of the ADMM for solving (21) in the following

$$\begin{cases} y^{k+1} \approx \arg\min_y \mathcal{L}_\sigma(y, \xi^k; x^k) + \frac{1}{2}\|y - y^k\|_S^2, \\ \xi^{k+1} \approx \arg\min_\xi \mathcal{L}_\sigma(y^{k+1}, \xi; x^k) + \frac{1}{2}\|\xi - \xi^k\|_T^2, \\ x^{k+1} = x^k - \tau\sigma\left(A^\top y^{k+1} + \xi^{k+1}\right), \end{cases} \tag{22}$$

where $\sigma > 0$ is a given penalty parameter, $\tau \in (0, \frac{1+\sqrt{5}}{2})$ is the dual steplength, which is typically chosen to be 1.618,

$$\mathcal{L}_\sigma(y, \xi; x) := \frac{1}{2}\|y\|^2 + \langle b, y \rangle + \kappa_\lambda^*(\xi) - \langle x, A^\top y + \xi \rangle + \frac{\sigma}{2}\|A^\top y + \xi\|^2$$

is the augmented Lagrangian function associated with problem (21) and $S$ and $T$ are two symmetric positive semidefinite matrices. The convergence results of such a general ADMM including the classical ones with the subproblems solved exactly have been discussed in Fazel et al. (2013) under some mild conditions. An inexact version for the general semi-proximal ADMM scheme and its convergence proof can be found in the recent paper by Chen et al. (2017).

In (22), the subproblem for updating $y$ can be handled by solving the following linear system corresponding to its optimality condition:

$$\left(I_m + \sigma AA^\top + S\right) y = Ax^k - b - \sigma A\xi^k + Sy^k.$$

The weight matrix $S$ in the proximal term can be simply chosen to be the zero matrix when $I_m + \sigma AA^\top$ admits a relatively cheap Cholesky factorization. Otherwise, we can adopt the rule elaborated in Subsection 7.1 in Chen et al. (2017) for choosing $S$ appropriately.

The subproblem for updating $\xi$ can be reformulated as:

$$\xi^{k+1} \approx \arg\min_\xi \kappa_\lambda^*(\xi)/\sigma + \frac{1}{2}\|\xi - (x^k/\sigma - A^\top y^{k+1})\|^2 + \frac{1}{2\sigma}\|\xi - \xi^k\|_T^2.$$

By utilizing the efficient algorithm for computing the proximal mapping $\text{Prox}_{\sigma\kappa_\lambda}(w^k)$ (Bogdan et al., 2015), together with the Moreau identity, we can choose $T = 0$ and update $\xi$ as follows:

$$\xi^{k+1} = \text{Prox}_{\kappa_\lambda^*/\sigma}(w^k/\sigma) = \frac{\left(w^k - \text{Prox}_{\sigma\kappa_\lambda}(w^k)\right)}{\sigma},$$

where $w^k := x^k - \sigma A^\top y^{k+1}$.

### 4.3 Results on solving the OSCAR model

In this subsection, we will test our proposed algorithm Newt-ALM for solving the OSCAR model and benchmark it against the APG algorithm implemented in the SLOPE package from Bogdan et al. (2015) and the ADMM scheme (Gabay and Mercier, 1976; Glowinski and Marrocco, 1975) presented in Subsection 4.2. The comparison among these three algorithms will be in terms of the computation time, the iteration number, and the accuracy measured via the relative KKT residual on several selected data from the UCI data repository Lichman (http://archive.ics.uci.edu/ml/datasets.html) and the BioNUS data set considered in Li et al. (2017). To demonstrate the performance of these three methods for solving the OSCAR model, with less consideration on the tuning parameter adjustment for pursuing a nice statistical behavior of the regularization model, here we manually choose the tuning parameters $w_1$ and $w_2$ as follows:

$$w_1 = a\|A^\top b\|_\infty \text{ and } w_2 = w_1/\sqrt{n} \tag{23}$$

with several testing values for the factor $a$. The sparsity is recorded in terms of the minimum number $k$ such that the first $k$ largest components in magnitude contribute a percentage of no less than $99.9\%$ for the $\ell_1$-norm. Results are shown in Table 1 and the data description is listed in Table 2. The "nnz" in Table 1 counts the number of nonzeros in the solution $x$ obtained by Newt-ALM such that $\text{nnz} = \min\{t : \sum_{i=1}^t |x_i| \geq 0.999\|x\|_1\}$

Table 1: The comparison results obtained by testing real data from the UCI and the BioNUS data sets with $(w_1, w_2)$ set as in (23). $A_1$: the ADMM with $\tau = 1.618$; $A_2$: the APG method implemented by the SLOPE package; $A_3$: our Newt-ALM.

| No. | a | nnz | $\eta$ $(A_1|A_2|A_3)$ | Time(s) $(A_1|A_2|A_3)$ | Iter No. $(A_1|A_2|A_3)$ |
|---|---|---|---|---|---|
| | 1e-7 | 1 | 4.37e-07\|8.47e-05\|5.26e-10 | 6.86\|3.01\|1.97 | 51\|25\|7 |
| p1 | 1e-8 | 2 | 2.53e-07\|1.56e-05\|3.37e-09 | 328.20\|3071.39\|7.11 | 2251\|32704\|37 |
| | 1e-9 | 158 | 9.72e-08\|7.37e-04\|1.83e-10 | 779.00\|4678.39\|25.33 | 4054\|50000\|52 |
| | 1e-7 | 1 | 4.00e-07\|4.53e-05\|2.60e-09 | 2.48\|1.78\|1.14 | 51\|25\|7 |
| p2 | 1e-8 | 3 | 5.88e-07\|3.14e-04\|3.61e-08 | 99.28\|961.80\|3.86 | 2101\|16525\|35 |
| | 1e-9 | 256 | 1.89e-07\|3.77e-04\|1.16e-11 | 277.84\|2605.04\|48.97 | 4627\|50000\|64 |
| | 1e-3 | 6 | 9.07e-09\|1.61e-07\|6.48e-09 | 549.90\|396.43\|4.07 | 7827\|5522\|23 |
| p3 | 1e-4 | 81 | 3.29e-08\|1.18e-06\|1.78e-09 | 4521.18\|3652.90\|8.77 | 50000\|50000\|39 |
| | 1e-5 | 100 | 3.87e-06\|3.26e-05\|2.62e-10 | 5781.13\|3776.61\|38.20 | 50000\|50000\|48 |
| | 1e-3 | 245 | 5.97e-10\|1.62e-07\|1.51e-11 | 351.95\|203.92\|9.22 | 802\|607\|7 |
| p4 | 1e-4 | 343 | 1.04e-07\|3.52e-07\|1.02e-09 | 2498.72\|7069.69\|21.94 | 4139\|21205\|23 |
| | 1e-5 | 419 | 1.84e-04\|3.69e-05\|3.63e-09 | 10800.80\|17138.70\|100.32 | 9893\|50000\|41 |
| | 1e-4 | 11 | 1.22e-07\|6.48e-07\|1.90e-08 | 144.99\|629.91\|2.69 | 1292\|17531\|21 |
| p5 | 1e-5 | 26 | 2.07e-07\|1.31e-05\|3.66e-09 | 815.35\|1618.06\|5.32 | 4896\|50000\|37 |
| | 1e-6 | 70 | 2.44e-07\|3.72e-04\|7.74e-09 | 859.92\|1613.81\|24.21 | 4320\|50000\|44 |
| | 1e-6 | 2 | 2.30e-07\|1.95e-08\|2.04e-10 | 215.07\|712.95\|4.31 | 1761\|10696\|23 |
| p6 | 1e-7 | 10 | 2.14e-07\|1.52e-05\|1.42e-09 | 583.87\|2899.70\|8.19 | 3509\|50000\|33 |
| | 1e-8 | 51 | 9.61e-09\|1.11e-04\|5.70e-10 | 1662.10\|2886.64\|19.88 | 6153\|50000\|46 |
| | 1e-3 | 8 | 2.41e-08\|1.70e-07\|2.04e-08 | 286.68\|89.86\|3.03 | 2713\|1401\|16 |
| p7 | 1e-4 | 39 | 9.98e-07\|5.40e-07\|4.96e-09 | 373.67\|1745.63\|4.64 | 1782\|29896\|24 |
| | 1e-5 | 120 | 9.75e-07\|6.81e-06\|1.18e-08 | 1303.22\|3032.56\|11.39 | 4266\|50000\|37 |
| | 1e-3 | 3 | 7.75e-08\|2.54e-07\|1.70e-07 | 2.20\|5.73\|0.73 | 302\|881\|13 |
| p8 | 1e-4 | 14 | 1.06e-07\|2.64e-06\|1.57e-07 | 4.95\|35.55\|0.95 | 513\|8013\|21 |
| | 1e-5 | 60 | 1.20e-07\|3.90e-06\|2.05e-08 | 16.28\|106.13\|1.61 | 1302\|31678\|32 |
| | 1e-3 | 1 | 1.75e-09\|1.42e-08\|2.56e-07 | 2.50\|1.01\|0.55 | 51\|28\|4 |
| p9 | 1e-4 | 6 | 6.64e-07\|3.99e-07\|5.24e-07 | 8.20\|22.47\|0.97 | 151\|715\|10 |
| | 1e-5 | 15 | 4.66e-07\|1.20e-06\|2.50e-08 | 46.97\|220.01\|2.02 | 701\|10905\|15 |
| | 1e-2 | 3 | 4.05e-08\|7.39e-07\|3.33e-08 | 4.20\|0.84\|0.72 | 623\|82\|20 |
| p10 | 1e-3 | 130 | 2.76e-08\|1.24e-06\|9.58e-09 | 13.94\|194.86\|3.13 | 1910\|25299\|33 |
| | 1e-4 | 160 | 9.66e-07\|2.75e-06\|1.53e-09 | 29.99\|220.75\|5.77 | 2963\|50000\|45 |
| | 1e-3 | 32 | 2.93e-08\|1.85e-06\|7.21e-09 | 42.61\|586.28\|2.41 | 3755\|50000\|34 |
| p11 | 1e-4 | 155 | 3.97e-07\|1.83e-05\|1.83e-08 | 32.15\|713.13\|3.25 | 2216\|50000\|40 |
| | 1e-5 | 193 | 2.71e-07\|4.08e-05\|8.81e-10 | 93.55\|352.37\|11.47 | 4608\|50000\|54 |
| | 1e-3 | 34 | 3.17e-08\|7.60e-07\|6.19e-09 | 9.84\|218.91\|1.06 | 3024\|42154\|34 |
| p12 | 1e-4 | 54 | 5.97e-07\|3.78e-06\|2.31e-09 | 14.02\|267.84\|2.11 | 3766\|50000\|39 |
| | 1e-5 | 96 | 7.55e-07\|7.58e-06\|8.24e-10 | 46.34\|267.88\|4.22 | 11501\|50000\|52 |
| | 1e-3 | 7 | 2.71e-08\|8.99e-06\|1.93e-08 | 524.93\|842.24\|2.17 | 38456\|50000\|44 |
| p13 | 1e-4 | 38 | 2.86e-08\|4.68e-05\|7.02e-09 | 761.06\|840.22\|5.83 | 50000\|50000\|45 |
| | 1e-5 | 151 | 2.28e-07\|1.33e-04\|3.05e-10 | 246.02\|670.25\|10.42 | 8130\|50000\|53 |
| | 1e-2 | 5 | 4.93e-08\|6.27e-07\|4.16e-08 | 2.20\|4.07\|0.62 | 997\|1187\|22 |
| p14 | 1e-3 | 33 | 4.65e-07\|6.93e-07\|5.91e-09 | 2.89\|68.43\|1.11 | 1302\|22498\|34 |
| | 1e-4 | 51 | 5.87e-07\|3.37e-06\|1.61e-08 | 7.89\|154.98\|2.08 | 4063\|50000\|39 |
| | 1e-2 | 1 | 5.46e-09\|4.47e-09\|5.08e-09 | 3.71\|1.72\|0.72 | 1002\|234\|23 |
| p15 | 1e-3 | 11 | 6.78e-08\|2.50e-06\|9.02e-09 | 9.89\|300.71\|1.17 | 2531\|50000\|32 |
| | 1e-4 | 76 | 1.40e-08\|1.79e-05\|7.98e-09 | 13.93\|304.82\|4.67 | 3144\|50000\|41 |

Table 2: The problem names and sizes

| No. | Problem name | $[m, n]$ |
|-----|--------------|----------|
| p1 | E2006.train | [16087,150360] |
| p2 | E2006.test | [3308,150358] |
| p3 | pyrim_scaled-expanded5 | [74, 201376] |
| p4 | triazines-scaled-expanded4 | [186,635376] |
| p5 | abalone_scale_expanded7 | [4177,6435] |
| p6 | bodyfat_scale_expanded7 | [252,116280] |
| p7 | housing_scale_expanded7 | [506,77520] |
| p8 | mpg_scale_expanded7 | [392, 3432] |
| p9 | space_ga_scale_expanded9 | [3107,5005] |
| p10 | DLBCL_H | [160, 7399] |
| p11 | lung_H1 | [203, 12600] |
| p12 | NervousSystem | [60, 7129] |
| p13 | ovarian_P | [253, 15153] |
| p14 | DLBCL_S | [47, 4026] |
| p15 | lung_M | [96, 7129] |

From Table 1, we observe that all the 45 tested instances are successfully solved by Newt-ALM within 2 minutes (for most of the cases within less than half a minute), while 3 and 22 cases have failed (i.e., not achieving our stopping criteria) to be solved by ADMM and SLOPE, respectively. Both the solution accuracy (as shown in the column under "$\eta$") and the computation time (as shown in the column under "Time(s)") show a tremendous computational advantage of Newt-ALM comparing to ADMM and SLOPE. In particular, for many of the instances corresponding to `p3, p4, p5, p6, p7, p13`, our algorithm can be more than 100 times faster than ADMM and SLOPE.

It is noteworthy that the dual-based ADMM also works better than SLOPE for a great majority of the tested instances. The performance profiles of these three algorithms for all 45 tested problems are presented in Figure 1. Recall that a point $(x, y)$ at a chosen algorithm curve implies that the algorithm can solve up to the desired accuracy $(100y)\%$ of all the tested instances within at most $x$ times of the fastest algorithm for each instance. More specifically, for $x = 150$, we can see from Figure 1 that even by consuming more than 150 times of the computation time taken by Newt-ALM, there are still around 40% and 10% of tested instances which are not successfully solved by SLOPE and ADMM, respectively.

### 4.4 The pathwise solution for a microarray data

The behavior of the OSCAR model for sparse feature selection and grouping for each specific instance relies heavily on the tuning parameters $w_1$ and $w_2$. To get a reliable and effective estimation for the coefficients of all involved predictors in the context of linear regression, a two-dimensional grid of various $w_1$ and $w_2$ values are tested to generate a solution path. The task of generating a solution path can be costly since each single pair of parameters $(w_1, w_2)$ will lead to a different instance of the OSCAR model. The path usually begins with appropriately chosen parameters that shrink all the coefficients to zero, and moves on until we are near the un-regularized solution by varying the values of the parameters. During the construction of the solution path, the warm start strategy (Friedman et al., 2007, 2010) is always used to accelerate the entire process by using the previous close-by solution as the initial point for the next problem.

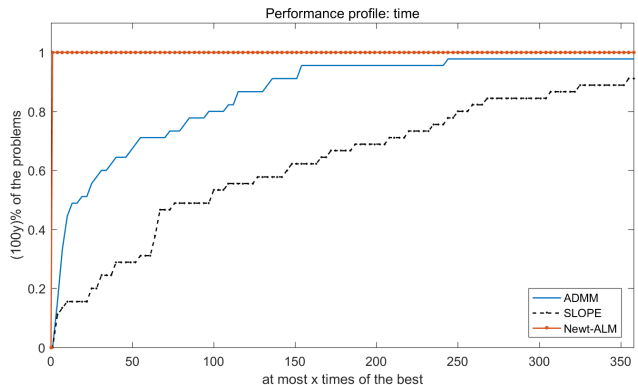

Figure 1: Time comparison for ADMM, SLOPE and Newt-ALM.

Table 3: Computation time comparison among Newt-ALM, ADMM and SLOPE for generating the partial solution paths, where the row "Ratio" reports the ratios of the computation time of each single algorithm to that of the fastest algorithm.

|         | Newt-ALM | ADMM   | SLOPE   |
|---------|----------|--------|---------|
| Time(s) | 27.74    | 323.65 | 1149.73 |
| Ratio   | 1        | 11.7   | 41.4    |

Here, we will use the microarray data set reported in Scheetz et al. (2006) and processed it by following Huang et al. (2008); Gu et al. (2017), where the design matrix $A \in \mathbb{R}^{m \times n}$ and the response vector $b \in \mathbb{R}^m$ with $m = 120$ and $n = 3000$. A partial solution path with the parameter $w_2$ of a fixed value $\|A^\top b\|_\infty / n^2$, and the parameter $w_1$ varying evenly in the interval $\left[10^{-4}, 10^{-2}\right] \times \|A^\top b\|_\infty$ for 100 different values will be constructed. The first 10 largest coefficients in magnitude of all the 100 numerical solutions are collected in Figure 2 by using ADMM, SLOPE and Newt-ALM, respectively. The timing comparison for generating the partial solution paths by these three algorithms is presented in Table 3.

As Figure 2 shows, all the three algorithms obtain almost the same partial solution paths for the microarray data in the above tuning parameter setting. This is due to fact that the size of the data is relatively small ($[m, n] = [120, 3000]$), and all the instances corresponding to the chosen tuning parameter pairs have rather sparse solutions (most of them have less than 10 nonzero components in the numerical solutions) and hence have been successfully solved by all three algorithms. Even for such a nice scenario, the Newt-ALM is still more than 10 times and 40 times faster than ADMM and SLOPE, respectively, as shown in Table 3. For difficult cases, such as large scale problems or those with relatively dense solutions in the high-dimensional linear regression, both the advantage in computation time and the estimation accuracy of our Newt-ALM will certainly be more significant, as can be observed in Table 1 and deduced from the computational complexity analysis in subsection 3.4 when the second-order sparsity of the generalized Jacobian has been fully uncovered and exploited.

With a two-dimensional grid of varying $w_1$ and $w_2$ values, we can also construct a three-dimensional scattergram to show those first $k$ (e.g., $k = 10$) largest components in magnitude for the microarray data. Figure 3 shows such a case, from which we can get a partial solution path along $w_1$ (or $w_2$) with any fixed $w_2$ (or $w_1$), or along any set of $(w_1, w_2)'s$ for the grid. Figure 3 shows the scattergram which collects the first 10 largest components in magnitude with $w_1$ and $w_2$ varying evenly in $\left[10^{-7}, 10^{-4}\right] \times \|A^\top b\|_\infty$ and $\left[10^{-4}, 10^{-2}\right] \times \|A^\top b\|_\infty/n$, respectively. All the $10,000$ problems are solved by our algorithm Newt-ALM in a total of about 70 minutes.

## 5. Discussions

In this paper we have proposed an efficient semismooth Newton-based augmented Lagrangian method for solving the OSCAR model in high-dimensional statistical regressions from the dual perspective. Numerical results have demonstrated the overwhelming superiority of the proposed algorithm on high-dimensional real data sets, comparing to the widely-used APG and ADMM. It is noteworthy that the original OSCAR model has been transformed to its dual counterpart before applying our method to take the advantage of the high-dimensional setting (i.e., the number of coefficients to be estimated is far larger than the sample size). The success of our second-order iterative method, both in accuracy and in computation time, relies heavily on the subtle second-order sparsity structure present in the generalized Jacobian matrix that corresponds to the second-order differential information of the underlying structured regularizer. Besides the least squares loss function adopted in the OSCAR model, our method is also applicable for the case of the logistic loss function, in which the desired nice properties of the corresponding subproblems are maintained to guarantee the efficiency and robustness of the algorithm. For classical statistical regression with larger sample size, our method is still applicable. But we may have to explore whether it is more efficient to apply our algorithmic framework directly to the OSCAR model, instead of our current application to the dual problem. The efficiency and effectiveness of our algorithm in solving high-dimensional linear regression with the OSCAR regularizer will greatly facilitate data analysis in statistical learning and related applications across a broad range of fields.

## Acknowledgments

We would like to acknowledge the financial support for the National Natural Science Foundation of China (11771038,11431002) and "111" Project of China (B16002), a start-up research grant from the Hong Kong Polytechnic University, and the Academic Research Fund (R-146-000-256-114) of the Ministry of Education of Singapore.

## Appendix A.

In this appendix we prove the following theorem from Section 2:

**Theorem** *Let* $\lambda \in \mathbb{R}_+^n$ *be such that* $\lambda = |\lambda|^\downarrow$. *Then* $\mathcal{M}(\cdot)$ *is a nonempty and compact valued, upper semicontinuous multifunction, and for any given* $y \in \mathbb{R}^n$, *every* $M \in \mathcal{M}(y)$ *is symmetric and positive semidefinite. Moreover, there exists a neighborhood* $U$ *of* $y$ *such that for all* $y' \in U$,

$$\mathrm{Prox}_{\kappa_\lambda}(y') - \mathrm{Prox}_{\kappa_\lambda}(y) - M(y' - y) = 0, \quad \forall M \in \mathcal{M}(y'). \tag{24}$$

∎

**Proof.** Let $y \in \mathbb{R}^n$ be an arbitrary point . Then it is obvious that $\mathcal{M}(y)$ is a nonempty and compact set. The symmetric and positive semidefiniteness of $M \in \mathcal{M}(y)$ is trivial by the definitions in (5) and (7). Now we claim that there exists a neighborhood $V$ of $y \in \mathbb{R}^n$ such that

$$\Pi^s(y') \subseteq \Pi^s(y), \quad \forall y' \in V.$$

This claim is trivial for $y = 0$ since $\Pi^s(0) = \mathbf{\Pi_n^s}$. For the case of a nonzero $y \in \mathbb{R}^n$, let $r$ be the number of distinct values in $|y|$, and $t_1$, ..., $t_r$ be all those distinct values satisfying $t_1 > t_2 > \cdots > t_r \geq 0$. Consider the following two cases:

Case I: If $t_r > 0$, set $\delta := \frac{1}{3} \min \left\{ t_r, \min_{1 \leq i \leq r-1} \{t_i - t_{i+1}\} \right\}$;

Case II: If $t_r = 0$, set $\delta := \frac{1}{3} \min_{1 \leq i \leq r-1} \{t_i - t_{i+1}\}$.

It is easy to verify that in both cases $\delta > 0$ and

$$\Pi^s(y') \subseteq \Pi^s(y), \quad \forall y' \in \mathbf{B}(y, \delta) \tag{25}$$

where $\mathbf{B}(y, \delta)$ is the 2-norm ball centered at $y$ with radius $\delta$. The upper semicontinuity of $\mathcal{M}$ then can be obtained from (25) and (6). The remaining part is to show (24). For any $y' \in \mathbf{B}(y, \delta)$ with $\delta$ defined as above, it is known from (4) and the inclusion property in (25) that

$$\mathrm{Prox}_{\kappa_\lambda}(y') - \mathrm{Prox}_{\kappa_\lambda}(y) = \pi^{-1} \left( x_\lambda(\pi y') - x_\lambda(\pi y) \right), \quad \forall \pi \in \Pi^s(y'). \tag{26}$$

By combining the properties in (6) and the fact that $\|\pi y' - \pi y\| = \|y' - y\|$, we know that there exists a neighborhood $U \subseteq \mathbf{B}(y, \delta)$ of $y$ such that for all $y' \in U$,

$$x_\lambda(\pi y') - x_\lambda(\pi y) = P(\pi y' - \pi y), \quad \forall P \in P(\pi y'), \forall \pi \in \Pi^s(y'),$$

which together with (26) leads to the desired result in (24). This completes the proof.

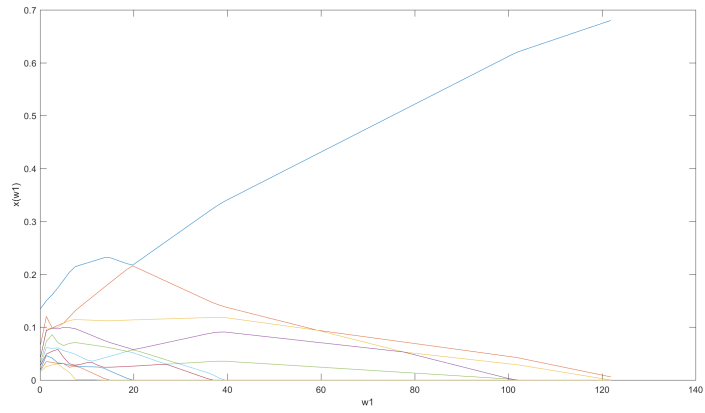

(a)

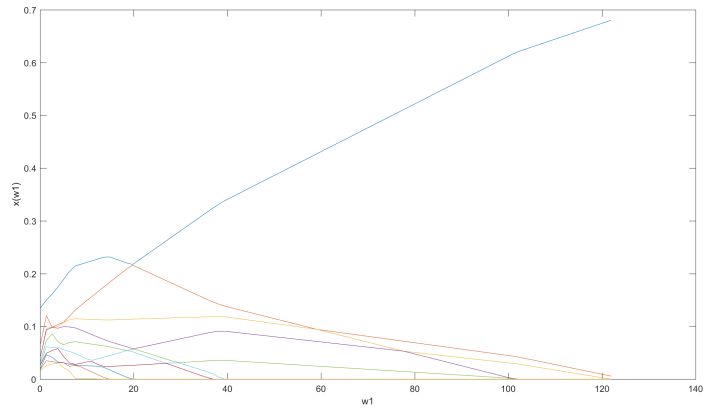

(b)

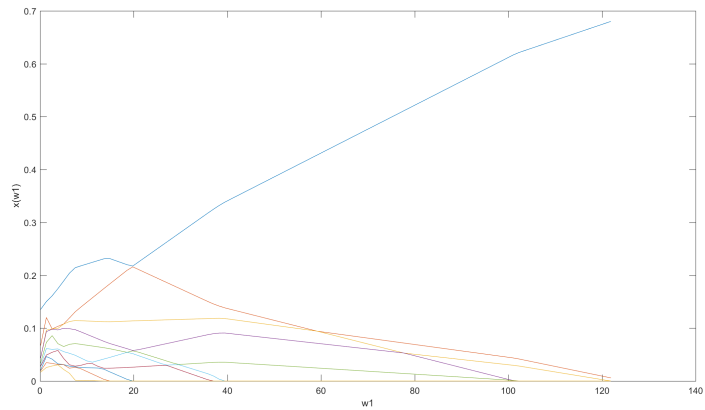

(c)

Figure 2: The partial solution paths with the first 10 largest coefficients in magnitude for the microarray data: (a) Newt-ALM; (b) ADMM; (c) SLOPE

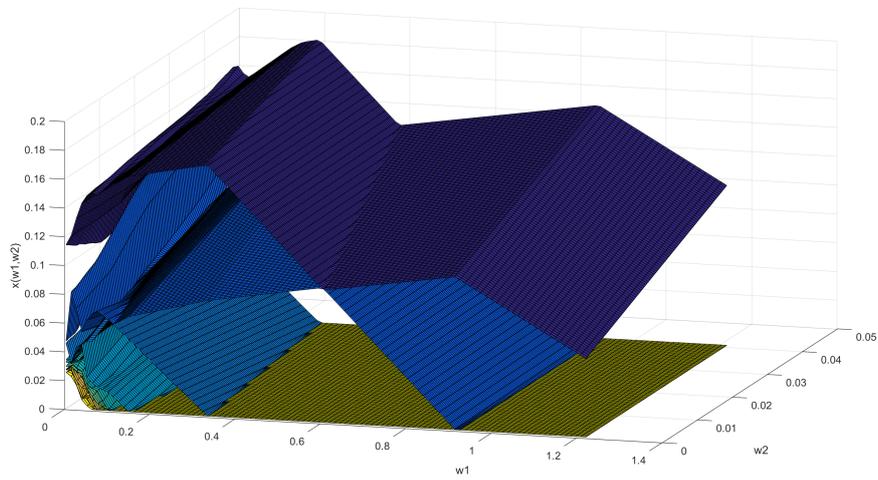

Figure 3: The first 10 largest components in magnitude of solutions with a two-dimensional grid of $w_1$ and $w_2$ values for the microarray data